\documentclass{ieeeconf}

\usepackage[backend=biber, style=ieee, isbn=false, doi=false]{biblatex}
\newcommand{\imgpath}{img}
\newcommand{\datapath}{img/data}

\usepackage{hyperref}
\hypersetup{
	colorlinks=true,
	linkcolor=black,
	citecolor=black,
	filecolor=black,
	urlcolor=black,
}

\DeclareSourcemap{\maps[datatype=bibtex]{
\map[overwrite]{
	\step[fieldsource=shortjournal,fieldtarget=journaltitle]
	\step[fieldsource=shortbooktitle,fieldtarget=booktitle]
}}}
\AtEveryBibitem{\clearfield{venue}}
\renewbibmacro*{event+venue+date}{}

\usepackage{basic_style}
\usepackage{commands}

\usepackage{balance}

\let\labelindent\relax
\usepackage[shortlabels]{enumitem}
\setlist[enumerate]{itemsep=0pt, itemindent=1em}

\usepackage[capitalise]{cleveref}
\crefname{section}{Sec.}{Secs.}
\Crefname{figure}{Fig.}{Figs.}

\newcommand*\circled[1]{\tikz[baseline=(char.base)]{
\node[shape=circle,draw,inner sep=1pt] (char) {#1};}}

\begin{document}	

\title{Low-Latency Digital Downconversion for Control Applications}

\author{Olof Troeng and Lawrence Doolittle%
\thanks{O. Troeng is with the Department of Automatic Control, Lund University, Sweden and is a member of the ELLIIT Strategic Research Area at Lund University. E-mail: \texttt{oloft@control.lth.se}.\qquad
L. Doolittle is with the Lawrence Berkeley National Lab, CA.
}%
}

\maketitle
\thispagestyle{empty}
	
\begin{abstract}
The slowly varying complex envelope of sinusoidal signals can be estimated in real-time using digital downconversion. In this paper,
we discuss the requirements on digital downconversion for control applications. Two low-latency downconversion filters are compared with respect to performance and implementation aspects.
\end{abstract}

\section{Introduction}
Estimating the complex envelope of digitized sinusoidal signals in real time is necessary in many control applications.
Some examples are: control of electromagnetic fields in particle accelerators \cite{Doolittle2006,Doolittle2008,Simrock2006,Schilcher2007,Hoffmann2008,Geng2019}, MEMS gyroscopes \cite{Saggin2019}, and laser stabilization \cite{Leibrandt2015}. 
For real-time computation of complex envelopes it is typical to use digital downconversion (DDC), see \cref{fig:down_conversion_complex_valued}.
There is a vast literature on DDC for telecommunications \cite{Crochiere1983,Vaidyanathan1993,Creaney2008,Ellingson2016}, but as we will see, it has limited applicability to control applications.

The focus of this paper is the filter $H(z)$ in Fig.~\ref{fig:down_conversion_complex_valued}. The filter $H(z)$ should both reject the double-frequency component of the mixer output and reduce aliasing. In telecommunications, its design is optimized with respect to passband flatness, stopband attenuation, and resource usage. For control applications, the latency, or more precisely the phase-drop at the feedback loop's cross-over frequency, is more important.


We present two types of low-latency DDC filters that are suitable for control applications: moving averaging and  two-sample reconstruction \cite{Doolittle2006,Doolittle2008}.
These filters have previously been used for field control in particle accelerators \cite{Doolittle2006,Doolittle2008,Simrock2006,Schilcher2007,Hoffmann2008}.
In this paper, we study their frequency-domain characteristics, compare their performance, and consider implementation aspects.
A key to our analysis is the realization that two-sampling reconstruction corresponds to \cref{fig:down_conversion_complex_valued} with a two-tap complex-coefficient filter $H(z)$.
The complex-valued perspective is particularly insightful in conjunction with complex-signal control analysis \cite{Cerezo2016,Troeng2017a}.


\scalebox{0.97}{\emph{Notation:} \mbox{The complex conjugate of $\bm{a}\!\in\! \mathbb{C}$ is denoted $\bm{a}^*$\!\!.}}

\emph{Remark:} Another interesting example of DDC, with its own set of challenges, is single-phase phasor estimation in three-phase power systems \cite{Phadke2017}.

\section{Background}

\subsection{Introduction to Digital Downconversion}

\begin{figure}
	\centering
	\tikzsetnextfilename{\imgpath/downconversion_complex_valued}%
	\input{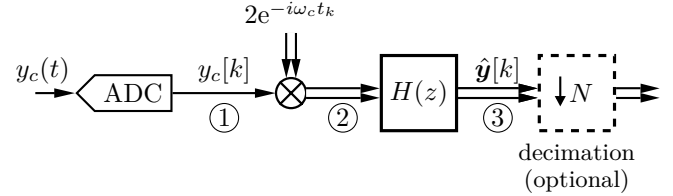}%
	
	\caption{Schematic of digital downconversion.
	The output $\hat{\by}$  is an estimate of the complex envelope of $y_c$.
	The filter $H(z)$, which is usually a low-pass filter, should reject the double-frequency components of the mixer output and act as an anti-aliasing filter (if the signal is decimated).
	Double arrows (%
	\protect\tikz[baseline=-2.5]{\protect\clip (-0.02,-0.2) rectangle (0.5, 0.2);  \protect\draw[-cplx] (0,0) -- (0.5,0);}) 
	indicate complex signals.	
	}
	\label{fig:down_conversion_complex_valued}
\end{figure}

A sinusoid with carrier frequency $\omega_c$, whose amplitude $A(t)$ and phase $\phi(t)$ vary slowly, has the form
\begin{multline}
y_c(t) = A(t) \cos(\omega_c t + \phi(t)) \\
= \Re\{ \underbrace{A(t) \me^{i\phi(t)}}_{=: \by(t)} \me^{i\omega_c t} \}  = \Re\{ \by(t) \me^{i\omega_c t} \}.
\end{multline}
The complex signal $\by(t) = A(t) \me^{i\phi(t)}$ is known as the \emph{(slowly varying) complex envelope} or the \emph{equivalent baseband signal}, of $y_c$. If $y_c(t)$ is sampled with a period $h$, at time instances $t_k = kh$, then a sampled version of $\by$ can be recovered by DDC as shown in Fig.~\ref{fig:down_conversion_complex_valued}.
\indent To better understand \cref{fig:down_conversion_complex_valued}, note that the signals at the indicated points are given by
\begin{itemize}
	\item[\circled{1}]
	\hspace*{0.5cm}$y_c[k] = \Re\!\left\{ \by[k] \me^{i\omega_c t_k} \right\}$\\[0.3em]
	\hspace*{1.31cm}$= (\by[k]\me^{i\omega_c t_k} + \conjnum{\by[k]}\me^{-i\omega_c t_k})/2$,
	\\[-0.5em]
	\item[\circled{2}] \hspace*{0.5cm} $\by[k] + \conjnum{\by[k]}\me^{-2i\omega_c t_k}$.
\end{itemize}
The low-pass filter removes the double-frequency component at $-2\omega_c$, giving the estimate $\hat{\by}[k] \approx \by[k]$ at \circled{3}.

DDC can be interpreted in the frequency domain as a translation by $-\omega_c$ followed by truncation of high frequencies, see \cite[Fig. 18.14]{Ellingson2016}.

\subsection{Baseband Model of DDC}
\begin{figure}
	\centering
	\tikzsetnextfilename{\imgpath/downconversion_baseband}%
	\input{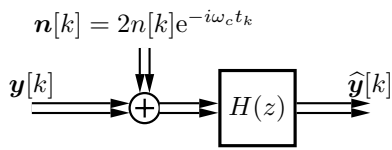}%

	\caption{
		Baseband model of DDC (cf. \cref{fig:down_conversion_complex_valued}).
		The double-frequency component  $\conjnum{\by}[k] \me^{-2i\omega_c t_k}$ is not included since it should be rejected by properly designed filters $H(z)$.
	}
	\label{fig:down_conversion_baseband}
\end{figure}

For control analysis it is convenient to transform the model in \cref{fig:down_conversion_complex_valued} to the  baseband. This yields the model in \cref{fig:down_conversion_baseband} that relates the complex envelope $\by$ of the input signal to the downconverted signal $\hat{\by}$.
In the baseband, it is seen that DDC amounts to adding complex-valued, cyclostationary noise $\bn$ and filtering by a linear, time-invariant (LTI) filter $H(z)$.

For control applications, the filter $H(z)$ should have low latency and hence  little impact on closed-loop performance.
It is nevertheless prudent to include the dynamics $H(z)$ and the measurement noise $\bn$ in analysis.


\subsection{DDC for Telecommunications}
\label{sec:ddc_communications_reqs}

In telecommunications, the signals of interest have a relatively flat spectrum over the channel bandwidth.
To estimate such signals without distortion, while rejecting spectral content outside the channel, calls for a filter $H(z)$ that is almost unity across the channel and has a fast transition to a well-attenuated stop band.
The desired features $H(z)$ can be summarized as
\begin{enumerate}[T1)]
	\item flat amplitude response in the passband, and
	\item excellent stopband suppression. \label{itm:stopband_suppr}
\end{enumerate}
In particular the second item requires significant engineering efforts since it is often difficult to achieve sufficient attenuation of adjacent channels and spurious interference\footnote{\label{fn:bob}Personal communication with B. Bernhardsson, former filter expert at Ericsson.}.

Typically, $H(z)$ is taken as an FIR filter $H(z) = h_0 + h_1 z^{-1} + \cdots + h_{N-1} z^{-(N-1)}$  with (conjugate-)symmetric
coefficients ($h_{N-1-k} = h_k^*$). The symmetry ensures that $H(z)$ is a \emph{linear-phase}
filter with the same phase response as a delay of $(N-1)/2$ samples.

The impulse and frequency responses of a typical filter $H(z)$ for telecommunications are shown in \cref{fig:impulse_response_wireless,fig:bodemag_wireless}.
\begin{figure}
	\centering
	\providecommand{\datafile}{\datapath/H_wireless_coeffs.csv}

	\tikzsetnextfilename{\imgpath/impulse_response_wireless}%
	\input{\imgpath/impulse_response_plot.tex}%

	\caption{
		Impulse response of a typical DDC filter $H(z)$ for telecommunications.
	}
	\label{fig:impulse_response_wireless}
	
	\vspace{1em}
	\renewcommand{\datafile}{\datapath/H_wireless_fr.csv}
	\newcommand{\xticks}{0, 50}

	\tikzsetnextfilename{\imgpath/bodemag_wireless}%
	\input{\imgpath/bodemag_filter.tex}%

	\caption{Magnitude response of the filter in \cref{fig:impulse_response_wireless}.
		The phase response equals that of a 42-sample time delay.}
	\label{fig:bodemag_wireless}
\end{figure}

In addition to the requirements T1 and T2, the design of $H(z)$ should allow an implementation with
\begin{enumerate}[T1), resume]
	\item low power consumption, and
	\item few hardware resources.
\end{enumerate}
For these reasons, the filter $H(z)$ is typically implemented as a cascade of filters at different rates, where the first one is typically a cascaded-integrator--comb filter \cite{Hogenauer1981}.
For details, see \cite{Creaney2008,Ellingson2016}.

To reduce the computations in the downstream baseband processing, the output of $H(z)$ is typically decimated as much as possible.
This can be done without loss of information due to that the spectrum of the baseband signal is centered around zero frequency.

In telecommunications, delay requirements are typically a minor concern$^{\ref{fn:bob}}$.

\section{DDC for Control Applications}
\label{sec:cavity_field_control}

\subsection{Requirements on DDC for Control Applications}

The requirements on DDC for control applications are quite different from those on DDC for telecommunications.
Typically, the signal of interest is very narrowband due to feedback, and hence it gives rise to ADC and mixer%
\footnote{If a heterodyne architecture is used.} harmonics, as well as a strong double-frequency component.
It is important to avoid that decimation folds these products onto the signal of interest. However, unlike in telecommunications, there is no need to worry about folding of spectral content from adjacent channels.

A main concern is instead to attenuate measurement noise (mostly ADC quantization noise) since it drives control errors and control signal activity. However, perhaps most important is to maintain sufficient phase margins.
The desired features of $H(z)$ for low-latency control applications could be summarized as
\begin{enumerate}[C1)]
	\item \scalebox{0.98}{small phase drop around the closed-loop bandwidth,}
	\item suppression of (ADC) measurement noise, \label{itm:noise_suppr}
	\item rejection of the double-frequency component,
	\item rejection of the DC-offset spur, and
	\item rejection of mixer and ADC harmonics\label{itm:harmonics}. 
\end{enumerate}

The trade-off between C1 and C2 is studied in the control literature \cite[Sec.~7.4]{Astrom1997}.
The following subsections introduce filters that address points C3--C5, although additional low-pass filtering might be necessary to achieve C2.
The filters $H(z)$ that we consider are rather simple, so resource usage should be a minor concern.

\emph{Remark:}
Requirement C2 is usually understood in a 2-norm sense,
while requirement T2 in \cref{sec:ddc_communications_reqs} is understood in a sup-norm sense.

\subsection{Decimation}
Since power and hardware resources tend to be less scarce in control applications, it could be feasible to use little, or no, decimation after the filter $H(z)$.
Decimation by a factor $N$ gives a controller period $h_\text{reg} = Nh$, which effectively corresponds to a delay $h_\text{reg}/2$ relative to continuous-time control \cite{Astrom1997}.
For this reason $h_\text{reg}$ should be kept small.
However, a small $h_\text{reg}$, increases the required resolution of coefficients and computations in the controller.
Trade-offs in the selection of $h_\text{reg}$ are  discussed in \cite{Astrom1997}.

\subsection{IQ Sampling}
\label{sec:iq_sampling}
A simple, low-latency approach to DDC is to take $f_c = f_s/4$ and \mbox{$H(z) = 1 + z^{-1}$}, where $f_c \coloneqq \omega_c/(2\pi)$ and $f_s \coloneqq 1/h$ is the sampling frequency.
This is approach is called \emph{IQ (in-phase and quadrature) sampling} \cite{Doolittle2006} or \emph{$f_s/4$ sampling} \cite{Ellingson2016}.
The zero of $H(z)$ in $-1$ eliminates the double-frequency component at $-f_s/2$.
IQ sampling can be implemented without multiplications and additions, which makes it attractive in terms of resource usage and power consumption \cite{Ellingson2016}.

A major drawback of IQ sampling is that odd mixer and ADC harmonics of the spectrally narrow input signal alias to the zero baseband frequency \cite{Doolittle2006}. This problem makes IQ sampling unsuitable for high-precision control.

\subsection{Non-IQ Sampling}
The aliasing problems of IQ sampling are avoided by using \emph{non-IQ (near-IQ) sampling} \cite{Doolittle2006,Schilcher2007} where $N$ samples are taken over $M$ periods of the carrier, i.e.,  $M/f_c = N/f_s$.

A filter $H(z)$ for non-IQ sampling should at least reject the double-frequency component at $-2\omega_c$.
Two such filter are presented in the next two subsections.

It will be convenient to introduce the phase advance between consecutive samples as
\begin{equation}
\Delta \coloneqq \omega_c h = 2\pi M/N.
\end{equation}
Note that $\Delta$ corresponds to $\omega_c$ in normalized angular frequency.

\subsection{Non-IQ DDC Filter: Moving Averaging}
\label{sec:moving_averaging}
For non-IQ sampling, an $N$-sample moving average filter
\begin{equation}
H_\text{\scriptsize MA}(z) \! = \! \frac{1}{N}\left( 1 + z^{-1} + \ldots + z^{-(N-1)} \right) \!=\! \frac{1 - z^{-N}}{1 - z^{-1}}
\label{eq:moving_average_filter}
\end{equation}
is a common choice \cite{Schilcher2007,Hoffmann2008}.

The impulse and frequency responses of an 11-sample moving average filter are shown in \cref{fig:impulse_response_CIC,fig:bodemag_cic}.
\begin{figure}
	\centering
	{
		\providecommand{\datafile}{\datapath/H_MA_11_coeffs.csv}
		\providecommand{\xmin}{0}
		\providecommand{\xmax}{26}
		\providecommand{\ymin}{0}
		\providecommand{\ymax}{0.14}

	\tikzsetnextfilename{\imgpath/impulse_response_MA_11}%
	\input{\imgpath/impulse_response_plot.tex}%

		\vspace{-0.3em}
		\caption{Impulse response of an 11-sample moving average filter.}
		\label{fig:impulse_response_CIC}
	}
	
	\vspace{1em}
	\centering
	\newcommand{\datafile}{\datapath/H_MA_11_fr.csv}

	\tikzsetnextfilename{\imgpath/bodemag_H_MA_11}%
	\input{\imgpath/bodemag_filter.tex}%

	\vspace{-0.3em}
	\caption{Magniude response of an 11-sample moving average filter; $\omega_s = 2\pi/h$.
		The phase response is that of a time delay of $11h/2$.}
	\label{fig:bodemag_cic}
	\vspace*{-0.7em}
\end{figure}
An attractive feature of the moving average filter \eqref{eq:moving_average_filter} is that it has zeros at the frequencies of the double-frequency component, the DC-offset spur, and all ADC and mixer harmonics (except multiples of $N$).

Decimation by a factor $N$ after the moving average filter  \eqref{eq:moving_average_filter} corresponds to a single-stage cascaded-integrator--comb (CIC) filter \cite{Hogenauer1981}.
The high side lobes in \cref{fig:bodemag_cic} indicate that aliasing could be a problem, but for decimation by $N$, these side lobes are neatly folded away from the zero frequency, enabling them to be rejected by a lower-rate low-pass filter 
\cite{Hogenauer1981}.

DDC with a moving average filter has two somewhat illuminating interpretations: as a short-time discrete Fourier transform with a single bin at frequency $\omega_c$ \cite{Hoffmann2008}; and as the least-squares estimate of $\by$ given $N$ observations $y_c[k], \ldots, y_c[k+N-1]$ \cite{Schilcher2007}.

\subsection{Non-IQ DDC Filter: Two-Sample Reconstruction}
\label{sec:two_sample}
\begin{figure}
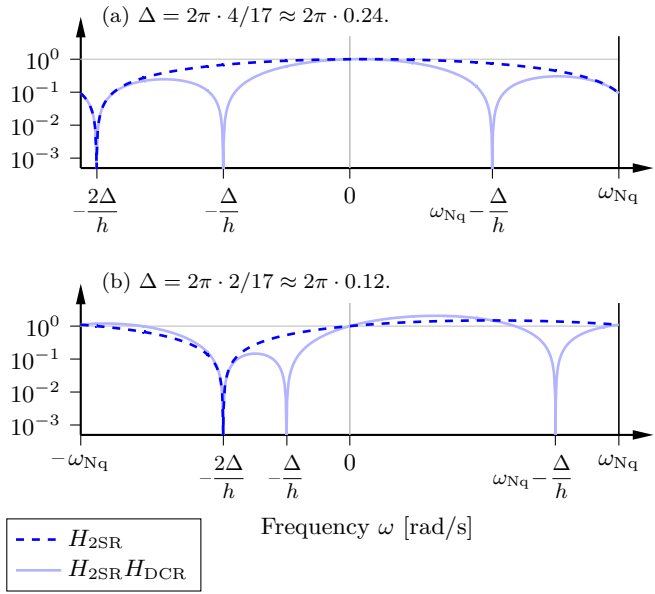

	\vspace{-0.3em}
	\centering
	\begin{minipage}{\columnwidth}
		\centering
		\providecommand{\xtick}{0, -2.35e7, -4.7e7, 2.65e7, 5e7}
		\newcommand{\datafile}{\datapath/H_2S_4_17_fr.csv}
		
		\newcommand{\xlabel}{}
		\providecommand{\captiony}{-0.75}		
		\newcommand{\captiontext}{(a) $\Delta = 2\pi \cdot 4/17  \approx 2\pi \cdot 0.24$.}				

	\tikzsetnextfilename{\imgpath/bodemag_2S_4_17}%
	\input{\imgpath/bodemag_filter_2S_dlin.tex}%

	\end{minipage}

	\vspace*{0.2em}		
	\begin{minipage}{\columnwidth}
		\centering

		\providecommand{\xtick}{0, -1.17e7, -2.35e7, 3.82e7, 5e7, -5e7}
		\newcommand{\datafile}{\datapath/H_2S_2_17_fr.csv}
		
		\newcommand{\showlegend}{}		
		\newcommand{\captiontext}{(b) $\Delta = 2\pi \cdot 2/17  \approx 2\pi \cdot 0.12$.}		

	\tikzsetnextfilename{\imgpath/bodemag_2S_2_17}%
	\input{\imgpath/bodemag_filter_2S_dlin.tex}%

	\end{minipage}
	\caption{Magnitude responses for the two-sample reconstruction filter $H_\text{2SR}(z)$ in \eqref{eq:2S_reconstruct} for two non-IQ ratios $\Delta$ (dashed lines). Solid lines show the inclusion of a (passband) DC rejection filter.}
	\label{fig:bodemag_two_sample}
\end{figure}

Two-sample reconstruction estimates the complex envelope $\by$ of a signal $y_c$ based on only two samples.
We start by showing that this method is of the form in \cref{fig:down_conversion_complex_valued} which allows for easier analysis than in \cite{Geng2019,Simrock2006,Doolittle2008}.
Assume that $y_c$ is sampled at times $t_{k-1} = (k-1)h$ and $t_k = kh$, and that $\by$ does not change between these two samples, i.e.,
\begin{subequations}
	\begin{equation} 
	\hspace*{-0.4em} y_c[k\!-\!1] \!=
	\Re\{\by \me^{i\omega_c t_{k-1}} \!\} \!
	= \! \frac{1}{2} \! \left( \by \me^{i\omega_c t_{k-1}} + \conjnum{\by}  \me^{-i\omega_c t_{k-1}}  \!\right)
	\label{eq:2S_kminusone}
	\end{equation}
	\begin{equation}
	\hspace*{-0.8em} y_c[k]  =
	\, \Re\{\by \me^{i\omega_c t_k} \} \,
	= \frac{1}{2}  \left( \by \me^{i\omega_c t_k} + \conjnum{\by}  \me^{-i\omega_c t_k}  \right).
	\label{eq:2S_k}
	\end{equation}
\end{subequations}
Taking \eqref{eq:2S_k} times $\me^{-i\omega_c t_k}$ minus \eqref{eq:2S_kminusone} times $\me^{-2i\Delta}\me^{-i\omega_c t_{k-1}}$ gives
\[
\me^{-i\omega_c t_{k}} y_c[k] - \me^{-2i\Delta} \me^{-i\omega_c t_{k-1}} y_c[k-1]  = \frac{1}{2} \left( \by - \me^{-2i\Delta} \by \right).
\]
Thus, $\by$ can be recovered as
\begin{equation}
\by = \frac{2}{1 - \me^{-2i\Delta}} \left( \me^{-i\omega_c t_{k}} y_c[k] - \me^{-2i\Delta} \me^{-i\omega_c t_{k-1}} y_c[k-1]  \right)
\end{equation}
which corresponds to filtering the signal $2\me^{-i\omega_c t_{k}} y[k]$ through the filter
\begin{align}
H_\text{2SR}(z) &=
\frac{1}{1 - \me^{-2i\Delta}} \left[1 - \me^{-2i\Delta} z^{-1} \right] \notag \\
&=
\frac{\me^{i\Delta}}{2i\sin{\Delta}} \left[1 - \me^{-2i\Delta} z^{-1} \right]
\label{eq:2S_reconstruct}
\end{align}
according to \cref{fig:down_conversion_complex_valued}.
The phase factor $\me^{i\Delta}/i$ may be discarded if only relative phase is of interest, but for our analysis it is convenient to have $\angle H_\text{2SR}(1) = 0$.
From the magnitude plot in \cref{fig:bodemag_two_sample} we see that $H_\text{2SR}(z)$ has a notch at frequency $-2\Delta/h = -2\omega_c$ (from the zero in $\me^{-2i\Delta}$)
that rejects the double-frequency component from the mixer.

Note from \eqref{eq:2S_reconstruct} that $H_\text{2SR}(z)$ has complex coefficients which implies that: (1) its frequency-response is not conjugate symmetric (see Fig~\ref{fig:bodemag_two_sample}) and (2) the complex-signal representation in \cref{fig:down_conversion_complex_valued} must be used to describe two-sample reconstruction, the version with two real-valued branches in, e.g., \cite[Sec.~18.4]{Ellingson2016}, does not apply.

From \eqref{eq:2S_reconstruct} we see that $\abs{\sin \Delta} \approx 1$ is necessary to avoid amplification of measurement noise \cite{Doolittle2008}. This corresponds to $M/N \approx 1/4$ (or $\approx 3/4$, etc), compare 
\makeatletter\cref@figure@name@plural\makeatother{} 
\ref{fig:bodemag_two_sample}a and \ref{fig:bodemag_two_sample}b.

For implementation details and Verilog code for two-sample reconstruction see \cite{Doolittle2008}.

\subsubsection*{DC-Spur Rejection}
Telecommunications-grade ADCs may have significant DC offset errors\footnote{\label{fn:comm-grade-adc}The offset error can be on the order of \SI{1}{\percent} of the maximum sine-wave amplitude.
For example,  LTC2175 from Linear Technologies has an offset error of up to 1.2\% of the maximum sine-wave amplitude.}.
From \cref{fig:down_conversion_baseband} it is seen that an offset error $n_0$ the generates a spur $2n_0 \me^{-i\omega_c t_k}$ in the baseband.
This spur can be removed by introducing a DC-rejection filter $(1 - z^{-2})/2$ before the digital mixer in \cref{fig:down_conversion_complex_valued} \cite{Doolittle2008}.
The DC-rejection filter is conveniently analyzed by transforming it to the baseband ($z \mapsto \me^{i\omega_c h} z = \me^{i\Delta} z$) where it corresponds to a filter $(1 - \me^{-2i\Delta} z^{-2})/2$.
Normalizing for unity gain at zero baseband frequency gives
\begin{equation}
H_\text{DCR}(z) = \frac{1 - \me^{-2i\Delta} z^{-2}}{1 - \me^{-2i\Delta}}
=
\frac{\me^{i\Delta}}{2i\sin{\Delta}} \left[1 - \me^{-2i\Delta} z^{-2} \right].
\label{eq:2P_DCR_filter}
\end{equation}
This filter should only be used if $\Delta \approx \pi/2$ to avoid amplifying measurement noise. Magnitude responses of the filter $H_\text{DCR}(z)$ together with $H_\text{2SR}(z)$ are shown in \cref{fig:bodemag_two_sample}.

An alternative is to reject the DC component at the mixer input by a high-pass filter $(z - 1)/(z-p)$ where $p$ is a real number slightly smaller than one. In the baseband, this filter corresponds to an IIR notch filter.
A filter of this type with $p=15/16$ is used for the field control systems at the Linac Coherent Light Source II.

\subsection{Non-IQ DDC Filter: IIR Notch Filter}
An alternative to moving averaging and two-sample reconstruction is to use a first-order IIR notch filter with a notch at $-2\omega_c$.
This approach avoids the latency of moving averaging and the constraint $\abs{\sin \Delta} \approx 1$ for two-sample reconstruction.
However, implementing an IIR filter at the sampling rate of the ADC could be technically challenging, and we are not aware that this approach has been used in practice. For this reason, we will not consider it further.

\section{Additional Low-Pass Filtering}
The noise rejection of the filters $H_*(z)$ in the previous section is often insufficient.
In this section we analyze two aspects related to additional low-pass filtering by a filter $F_\text{LP}(z)$, see \cref{fig:ddc_lp}.

\subsection{Moving Averaging vs. Two-Sample Reconstruction}
\begin{figure}
	\centering
	\tikzsetnextfilename{\imgpath/ddc_lp}%
	\input{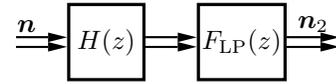}%

	\caption{A low-pass filter  $F_\text{LP}(z)$ (\cref{eq:def_F_LP}) after the DDC filter $H(z)$ improves the noise rejection from $\bn$ to  $\bn_2$.
	}
	\label{fig:ddc_lp}
\end{figure}

Which of the two filters in the previous section is the better choice? Two-sample reconstruction has a shorter latency, but moving averaging provides better high-frequency roll off. For a comparison, we considered the setup in \cref{fig:ddc_lp} where $\bn$ is white noise and  $F_\text{LP}(z)$ is a first-order low-pass filter with bandwidth $\omega_\text{LP}$,
\begin{equation}
F_\text{LP}(z) = \frac{1-a}{1 -az^{-1}}, \qquad a = \me^{-\omega_\text{LP} h}.
\label{eq:def_F_LP}
\end{equation}

\begin{figure}
	\centering
	\tikzsetnextfilename{\imgpath/bode_filter_comparison}%
	\input{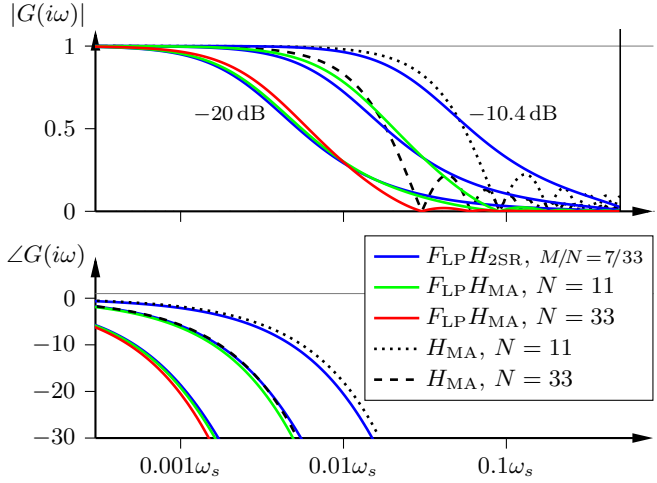}%

	\caption{
		Frequency responses for three groups of filter combinations $F_\text{LP}H_\ast$ with the $\omega_\text{LP}$ tuned for noise rejections $\norminline{F_\text{LP}H_\ast}_2^2$ of \SI{-20}{dB}, \SI{-15.2}{dB}, and \SI{-10.4}{dB}. The dashed line corresponds to a 33-tap moving average filter (\SI{-15.2}{dB}) and the dotted line to an 11-tap moving average filter (\SI{-10.4}{dB}).  Approximate values of the bandwidths  $\omega_\text{LP}$ in the three groups were $0.00032\omega_s$, $0.01\omega_s$, and $0.032\omega_s$. The frequency response of $H_\text{2SR}$ is close enough to conjugate symmetric that we do not show negative frequencies.
	}
	\label{fig:bode_filter_comparison}
\end{figure}

For each filter $H_*(z)$, the lowpass bandwidth $\omega_\text{LP}$ was tuned for three levels of noise rejection $\norm{H_*(z)F_\text{LP}(z)}_2^2$.
The frequency responses of the resulting filter combinations are shown in \cref{fig:bode_filter_comparison}.
It is seen that two-sample reconstruction enables slightly faster roll off with less phase retardation, which enables better feedback performance. The difference is quite small at low frequencies, but becomes noticeable if a high cross-over frequency is desired. In particular, for long moving average filters (note the dashed line in \cref{fig:bode_filter_comparison}).

That moving averaging is usually suboptimal for control applications is well-known and can be intuitively understood from that equal weight is placed on the most recent sample and the $N$th most recent sample.

\subsection{Order of Low-Pass Filtering and Decimation}
If the output signal $\hat{\by}$ is to be decimated, one may ask whether this should be done before or after the low-pass filter $F_\text{LP}(z)$ (\cref{fig:filt_dec_order_illustration}).
Decimating by a factor $N$ before low-pass filtering enables a lower-rate filter implementation $\widetilde{F}_\text{LP}(z)$ (as in \eqref{eq:def_F_LP} but with $a = \me^{-\omega_\text{LP} N h}$), but leads to increased aliasing. However, due to the nice folding properties of moving averaging \eqref{eq:moving_average_filter} and two-sample reconstruction \eqref{eq:2S_reconstruct} (if $\abs{\sin \Delta} \approx 1$), this is not overly problematic.

\Cref{fig:filt_dec_order_plot} shows how the order of filtering and decimation affects the resulting noise level as a function of the bandwidth $\omega_\text{LP}$ of the low-pass filter. The output variance of the periodic system in \cref{fig:filt_dec_order_plot}b was computed as $\norminline{\widetilde{F}_\text{LP}(z^N) H(z)}_2^2$ (see \emph{noble identities} in \cite{Vaidyanathan1993}).

In the case of moving average filters, we see from \cref{fig:filt_dec_order_plot}a that the order of low-pass filtering and decimation has little impact on the noise rejection if the bandwidth $\omega_\text{LP}$ is lower than \SI{100}{kHz}.
For long moving average filters (and consequently high decimation ratios) it is seen that there is much to gain from low-pass filtering with $F_\text{LP}(z)$ before decimation.

For two-sample reconstruction, we see from \cref{fig:filt_dec_order_plot}b that the order of decimation and low-pass filtering makes little difference if $\abs{\sin \Delta} \approx 1$ (which was necessary for two-sample reconstruction to be viable in the first place, see \cref{sec:two_sample}).
\begin{figure}
	\centering
	\begin{minipage}[c]{0.49\columnwidth}
	\centering
	\scalebox{0.88}{
	\tikzsetnextfilename{\imgpath/ddc_lp_dec}%
	\input{\imgpath/ddc_lp_dec.tex}%
	
	}
	
	\vspace{1em}
	{\footnotesize (a) Decimation after $F_\text{LP}(z)$.}
\end{minipage}
\begin{minipage}[c]{0.49\columnwidth}
	\centering
	\scalebox{0.88}{
	\tikzsetnextfilename{\imgpath/ddc_dec_lp}%
	\input{\imgpath/ddc_dec_lp.tex}%

	}	
	
	\vspace{1em}
	{\footnotesize (b) Decimation before $F_\text{LP}(z)$.}		
\end{minipage}
	
	\caption{Two alternatives for implementing additional low-pass filtering.
	Approach (a) reduces aliasing but requires the filter $F_\text{LP}(z)$ to run at a higher rate than in (b). The coefficients of $F_\text{LP}(z)$ and $\widetilde{F}_\text{LP}(z)$ differ due to the different sampling rates.
	}
	\label{fig:filt_dec_order_illustration}
	\centering
	{

	\tikzsetnextfilename{\imgpath/filt_dec_order_norm_moving_average}%
	\input{\imgpath/filt_dec_order_norm.tex}%

	{\footnotesize (a) Moving averaging and decimation by $N$. }
	}
	
	{
	\newcommand{\twosampleplot}{}

	\tikzsetnextfilename{\imgpath/filt_dec_order_norm_2_sample}%
	\input{\imgpath/filt_dec_order_norm.tex}%

	{\footnotesize (b) Two-sample reconstruction and decimation by 2. }
	}
	\caption{Noise rejection from filtering by $H(z)$ and  $F_\text{LP}(z)$ (\cref{fig:filt_dec_order_illustration}) with decimation \emph{after} low-pass filtering (solid lines) and decimation \emph{before} low-pass filtering (dashed lines), relative to only low-pass filtering by $F_\text{LP}(z)$ at the original sample period $h$.
}
\label{fig:filt_dec_order_plot}

\end{figure}

\section{Examples}
\subsubsection*{LCLS-II} The field control systems for the Linac Coherent Light Source II use non-IQ sampling with $M/N = 7/33$ ($f_c = \SI{20}{MHz}$ and $f_s =   \SI{94.29}{MHz}$) together with two-sample reconstruction \cite{Doolittle2016}. Due to the high feedback gain ($K \approx 1000$), an additional first-order low-pass filter with a bandwidth between \SI{50}{kHz}--\SI{200}{kHz} is used to limit the control signal activity.

\subsubsection*{ESS} The field control systems for the European Spallation Source use non-IQ sampling with $M/N = 3/14$ ($f_c = \SI{25.16}{MHz}$ and $f_s =   \SI{117.40}{MHz}$) together with moving averaging and decimation by $N$.

\section{Conclusions}
We have discussed DDC for control applications and how it differs from DDC for telecommunications.
We considered two types of DDC filters that are suitable for control applications: moving averaging, which is a decent and trouble-free choice; and two-sample reconstruction which allows slightly better performance, in particular for extremely fast feedback. However, with two-sample reconstruction it is necessary to be mindful of harmonics, the DC-offset spur, and that $\abs{\sin \Delta} \approx 1$.

Typically, the output of these two filters needs to be low-pass filtered and decimated. It is preferable to low-pass filter before decimating, but the noise increase from doing it the other way around is small, except for long moving averaging filters or if little filtering is applied.

\section{Acknowledgment}
Bo Bernhardsson is thanked for sharing his expertise on filter design for telecommunications and for his helpful comments and suggestions.

\printbibliography

\end{document}